\documentclass[runningheads]{llncs}
\usepackage{amsmath,amssymb,mathtools}
\usepackage[backgroundcolor=black!20]{todonotes}
\usepackage{hyperref,cleveref}
\usepackage{algorithm}
\usepackage[noend]{algpseudocode}
\usepackage{nicefrac}
\usepackage{float}
\usepackage{booktabs}
\usepackage[inline]{enumitem}
\usepackage{graphicx}
\usepackage{caption}
\usepackage{grffile} 
\usepackage{placeins}

\usepackage{color}
\usepackage{setspace}

\algrenewcommand\algorithmicrequire{\textbf{Input:}}
\algrenewcommand\algorithmicensure{\textbf{Output:}}

\newcommand{\eps}{\varepsilon}
\newcommand{\conv}{\operatorname{conv}}
\newcommand{\R}{\mathbb{R}}

\newcommand{\sprod}[2]{\langle #1, #2 \rangle}
\newcommand{\zero}{\mathbf{0}}
\newcommand{\one}{\mathbf{1}}
\newcommand{\opt}{\mathrm{OPT}}
\newcommand{\ub}{\mathrm{UB}}

\newcommand{\oracle}{\mathcal{A}}

\newcommand{\code}[1]{\texttt{#1}}
\newcommand{\pot}{\Phi}

\makeatletter
\renewcommand{\@Opargbegintheorem}[4]{%
  #4\trivlist\item[\hskip\labelsep{#3#2\@thmcounterend}]}
\makeatother

\begin{document}

\title{A Simple Method for Convex Optimization in the Oracle Model} 

\author{Daniel Dadush\inst{1}\thanks{This project has received funding from the European Research Council (ERC) under the European Union's Horizon 2020 research and innovation programme (grant agreement QIP--805241)} \and 
Christopher Hojny\inst{2} \and 
Sophie Huiberts\inst{1} \and
Stefan Weltge\inst{3}} 

\authorrunning{Dadush, Hojny, Huiberts, Weltge}

\institute{Centrum Wiskunde \& Informatica, Netherlands\\
  \email{\{dadush,s.huiberts\}@cwi.nl}
  \and
  Eindhoven University of Technology, Netherlands\\
  \email{c.hojny@tue.nl}
  \and
  Technical University of Munich, Germany\\
  \email{weltge@tum.de}
}

\maketitle

\begin{abstract}
We give a simple and natural method for computing approximately optimal
solutions for minimizing a convex function $f$ over a convex set $K$
given by a separation oracle. Our method utilizes the Frank--Wolfe algorithm
over the cone of valid inequalities of $K$ and subgradients of $f$. Under the
assumption that $f$ is $L$-Lipschitz and that $K$ contains a ball of radius
$r$ and is contained inside the origin centered ball of radius $R$, using
$O(\frac{(RL)^2}{\eps^2} \cdot \frac{R^2}{r^2})$ iterations and calls to the
oracle, our main method outputs a point $x \in K$ satisfying $f(x) \leq \eps
+ \min_{z \in K} f(z)$.

Our algorithm is easy to implement, and we believe it can serve as a useful
alternative to existing cutting plane methods. As evidence towards this, we
show that it compares favorably in terms of iteration counts to the standard
LP based cutting plane method and the analytic center cutting plane method,
on a testbed of combinatorial, semidefinite and machine learning instances.  

\end{abstract}

\section{Introduction}

We consider the problem of minimizing a convex function $f\colon \R^n \to \R$ over a compact convex set $K \subseteq \R^n$.
We assume that $K$ contains an (unknown) Euclidean ball of radius $r > 0$ and is contained inside the origin centered ball of radius $R > 0$, and that $f$ is $L$-Lipschitz.
We have first-order access to $f$ that yields $f(x)$ and a subgradient of $f$ at $x$ for any given $x$.
Moreover, we only have access to $K$ through a separation oracle (SO), which, given a point $x \in \R^n$, either asserts that $x \in K$ or returns a linear constraint valid for $K$ but violated by $x$.

Convex optimization in the SO model is one of the fundamental settings in optimization.
The model is relevant for a wide variety of implicit optimization problems, where an explicit description of the defining inequalities for $K$ is either too large to store or not fully known.
The SO model was first introduced in~\cite{YN83} where it was shown that an additive $\eps$-approximate solution can be obtained using $O(n\log(LR/(\eps r)))$ queries via the center of gravity method and $O(n^2\log(LR/(\eps r)))$ queries via the ellipsoid method.
This latter result was used by Khachiyan~\cite{Khachiyan79} to give the first polynomial time method for
linear programming.
The study of oracle-type models was greatly extended in the classic book of Grötschel, Lovász, and Schrijver~\cite{GLS}, where many applications to combinatorial optimization were provided.
Further progress on the SO model was given by Vaidya~\cite{Vaidya96}, who showed that the $O(n \log(LR/(\eps r)))$ oracle complexity can be efficiently achieved using the so-called volumetric barrier as a potential function, where the best current running time for such methods was given very
recently~\cite{LSW15,JLSW20}. 

From the practical perspective, two of the most popular methods in the SO
model are the standard linear programming (LP) based cutting plane method,
independently discovered by Kelley~\cite{Kelley60},
Goldstein-Cheney~\cite{CheneyGoldstein59} as well as Gomory~\cite{Gomory58}
(in the integer programming context), and the analytic center cutting plane
method~\cite{Sonnevend88} (ACCPM). 

The LP based cutting plane method, which we henceforth dub the standard cut
loop, proceeds as follows: starting with finitely many linear underestimators
of $f$ and linear constraints valid for $K$, in each iteration it solves a
linear program that minimizes the lower envelope of $f$ subject to the
current linear relaxation of $K$. The resulting point $x$ is then used to
query $f$ and the SO to obtain a new underestimator for $f$ and a new
constraint valid for $K$. Note that if $f$ is a linear function, it
repeatedly minimizes $f$ over linear relaxations of $K$. While it is
typically fast in practice, it can be unstable, and no general quantitative
convergence guarantees are known for the standard cut loop.

To link to integer programming, in that context $K$ is the convex hull of
integer points of some polytope $P$ and the objective is often linear, and
the method is initialized with a linear description of $P$. A crucial
difference there is that the separator SO is generally only efficient
when queried at vertices of the current relaxation. 

ACCPM is a barrier based method, in which the next query point is the
minimizer of the barrier for the current inequalities in the system. ACCPM is
in general a more stable method with provable complexity guarantees.
Interestingly, while variants of ACCPM with $O(n \log(RL/(r\eps))^2)$ convergence
exist, achieved by judiciously dropping constraints~\cite{AtkinsonVaidya95},
the more practical variants have worse guarantees.
For instance, if $K$ is the ball of radius $R$, the standard variant of ACCPM is only shown to achieve $O(n(RL/\eps)^2 \log(RL/\eps))$ convergence~\cite{Nesterov95}.

In this paper, we describe a new method for convex optimization in the SO
model that computes an additive $\eps$-approximate solution within
$O(\nicefrac{R^4L^2}{r^2 \eps^2})$ iterations. Our algorithm is easy to
implement, and we believe it can serve as a useful alternative to existing
methods. In our experimental results, we show that it compares favorably in
terms of iteration counts to the standard cut loop and the analytic center
cutting plane method, on a testbed of combinatorial, semidefinite and machine
learning instances.

Before explaining our approach, we review the relevant work in related models.
To begin, there has been a tremendous amount of work in the context of
first-order methods~\cite{ben2001lectures,beck2017first}, where the goal is to minimize a possibly complicated
function, given by a gradient oracle, over a \emph{simple domain} $K$
(e.g., the simplex, cube, $\ell_2$ ball). These methods tend to have cheap iterations and to achieve ${\rm poly}(1/\eps)$ convergence rates. They are often superior in practice when the requisite accuracy is low or moderate, e.g., within $1\%$ of optimal.  For these methods, often variants
of (sub-)gradient descent, it is generally assumed that computing (Euclidean)
projections onto $K$ as well as linear optimization over $K$ are easy.
 If one only assumes access to
a linear optimization (LO) oracle on $K$, $K$ can become more interesting
(e.g., the shortest-path or spanning-tree polytope). In this
context, one of the most popular methods is the so-called
\emph{Frank--Wolfe} algorithm~\cite{FW56} (see~\cite{Jaggi13} for a modern
treatment), which iteratively computes a convex combination of vertices of
$K$ to obtain an approximate minimizer of a smooth convex function.

In the context of combinatorial optimization, there has been a considerable
line of work on solving (implicit) packing and covering problems using the
so-called multiplicative weights update (MWU)
framework~\cite{SM90,PST95,GK07}. In this framework, one must be able to
implement an MWU oracle, which in essence computes optimal solutions for the
target problem after the ``difficult'' constraints have been aggregated
according to the current weights. This framework has been applied for getting
fast $(1 \pm \eps)$-approximate solutions to multi-commodity
flow~\cite{SM90,GK07}, packing spanning trees~\cite{CQ17a}, the Held--Karp
approximation for TSP~\cite{CQ17b}, and more, where the MWU oracle
computes shortest paths, minimum cost spanning trees, minimum cuts
respectively in a sequence of weighted graphs. The MWU oracle is in general
just a special type of LO oracle, which can often be interpreted as a SO that returns a maximally violated
constraint. While certainly related to the SO model, it is not entirely clear
how to adapt MWU to work with a general SO, in particular in settings
unrelated to packing and covering.

A final line of work, which directly inspires our work, has examined simple
iterative methods for computing a point in the interior of a cone $\Sigma$
that directly apply in the SO model. The application of simple iterative
methods for solving conic feasibility problems can be traced to Von Neumann
in 1948 (see~\cite{Dantzig91}), and a variant of this method, the perceptron
algorithm~\cite{Rosenblatt58} is still very popular today. Von Neumann's
algorithm computes a convex combination of the defining inequalities of the
cone, scaled to be of unit length, of nearly minimal Euclidean norm. The
separation oracle is called to find an inequality violated by the current
convex combination, and this inequality is then used to make the current convex
combination shorter, in an analogous way to Frank--Wolfe. This method is guaranteed to
find a point in the cone in $O(1/\rho^2)$ iterations, where $\rho$ is the
so-called width of $\Sigma$ (the radius of the largest ball contained in
$\Sigma$ centered at a point of norm $1$). Starting in 2004, polynomial time
variants of this and related methods (i.e., achieving $\log (1/\rho)$
dependence) have been found~\cite{Betke04,DV08,Chubanov12}, which iteratively
``rescale'' the norm to speed up the convergence. These rescaled variants can
also be applied in the oracle setting~\cite{BFV09,Chubanov17,DVZ20} with
appropriate adaptations. The main shortcoming of existing conic approaches is that they
are currently not well-adapted for solving optimization problems rather than feasibility problems.

\paragraph{Our approach.}
In this work, we build upon von Neumann's approach and utilize the Frank--Wolfe algorithm
over the cone of valid inequalities of $K$ as well as the subgradients of $f$ in a way that yields a clean, simple, and flexible framework for solving general convex optimization problems in the SO model.
For simpler explanation, let us assume that $f(x) = \langle c,x \rangle$ is a linear function and that we know an upper bound $\ub$ on the minimum of $f$ over $K$.
Given some linear inequalities $\langle a_i,x \rangle \le b_i$ that are valid for all $x \in K$,
our goal is to find convex combinations $p$ of the \emph{homogenized} points $(c,\ub)$ and $(a_i,b_i)$ that are ``close'' to the origin.
Note that if $p = \zero$, the fact that $K$ is full-dimensional implies that $(c,\ub)$ appears with a nonzero coefficient and hence $(-c,-\ub)$ is a nonnegative combination of the points $(a_i,b_i)$, which in turn shows that $\ub$ is equal to the minimum of $f$ over $K$.
In view of this, we will consider a potential $\pot\colon \R^{n+1} \to \R_+$ with the property that if $\pot(p)$ is sufficiently small, then the convex combination will yield an explicit certificate that $\ub$ is close to the minimum of $f$ over $K$.

Given a certain convex combination $p$, note that the gradient of $\pot$ at $p$ provides information about whether moving towards one of the known points will (significantly) decrease $\pot(p)$.
However, if no such known point exists, it turns out that the ``dehomogenization'' of the gradient (a scaling of its projection onto the first $n$ coordinates) is a natural point $x \in \R^n$ to query the SO with.
In fact, if $x \in K$, it will have improved objective value with respect to $f$.
Otherwise, the SO will provide a linear inequality such that moving towards its homogenization decreases $\pot(p)$.

In this work, we will show that the above paradigm immediately yields a rigorous algorithm for various natural choices of $\pot$ and scalings of inequalities.
We will also see that general convex functions can be directly handled in the same manner by simply replacing $(c,\ub)$ with all subgradient cuts of $f$ learned throughout the iterations.
The same applies to pure feasibility problems for which we set $f = \zero$.
The convergence analysis of our algorithm is simple and based on standard estimates for the Frank--Wolfe algorithm.

Besides its conceptual simplicity and distinction to existing methods for convex optimization in the SO model, we also regard it as a practical alternative.
In fact, in terms of iterations, our vanilla implementation in \code{Julia}\footnote{\url{https://github.com/christopherhojny/supplement_simple-iterative-methods-linopt-convex-sets}} performs similarly and often even better than the standard cut loop and the analytic center cutting plane method evaluated on a testbed of oracle-based linear optimization problems for matching problems, semidefinite relaxations of the maximum cut problem, and LPBoost.
Moreover, the flexibility of our framework leaves several degrees of freedom to obtain optimized implementations that outperform our naive implementation.

\paragraph{Acknowledgments}
We would like to thank Robert Luce and Sebastian Pokutta for their very valuable feedback on our work.

\section{Algorithm}
\label{sec:alg}

Recall that we are given first-order access to a convex function $f\colon \R^n \to \R$ that we want to minimize over a convex body $K \subseteq \R^n$.
In the case where $f$ is not differentiable, with a slight abuse of notation we interpret $\nabla f(x)$ to be \emph{any} subgradient of $f$ at $x$.
We can access $K$ by a separation oracle that, given a point $x \in \R^n$, either asserts that $x \in K$ or returns a point $(a,b) \in \oracle \subseteq \R^{n+1}$ with $\langle a,x \rangle > b$ such that $\langle a,y \rangle \le b$ holds for all $y \in K$.
Here, $\langle \cdot,\cdot \rangle$ denotes the standard scalar product and we assume that all points in $\oracle$ correspond to linear constraints valid for $K$.
To state our algorithm, let $\|\cdot\|$ denote any norm on $\R^{n+1}$ and $\|\cdot\|_*$ its dual norm.
Moreover, let $\pot\colon \R^{n+1} \to \R_+$ be any strictly convex and differentiable function
with $\min_{x \in \R^{n+1}} \pot(x) = \pot(0) = 0$.
Our method is given in~\Cref{algorithm}, in which we denote the number of iterations by $T$ for later reference.
However, $T$ does not need to be specified in advance, and the algorithm may be stopped at any time, e.g., when a solution or bound of desired accuracy has been found.

\begin{algorithm}[!h]
    \caption{}
    \label{algorithm}
    \begin{algorithmic}[1]
        \State $\ub \gets \infty$,
        $A_1 \gets \{(\zero,1)/\|(\zero,1)\|_*\}$,
        $G_1 \gets \emptyset$
        \For{$t = 1,2,\dots,T$}
        \State \label{algmin}$p_t \gets \arg\min\{ \pot(p) : p \in \conv(A_t \cup G_t)\}$
        \If{$p_t = \zero$}
            \Return $\ub$.
        \EndIf
        \State \label{algdefx}$x_t \gets - \nabla \pot(p_t)[1:n]/ \nabla \pot(p_t)[n+1]$
        \If{$x_t \in K$}
            \State $\ub \gets \min \{ \ub, f(x_t) \}$
            \State $A_{t+1} \gets A_t$.
            \State $G_{t+1} \gets G_t \cup \{(\nabla f(x_t), \sprod{\nabla f(x_t)}{x_t})\}$
        \Else
            \State get $(a,b) \in \oracle$, with $\sprod a {x_t} > b$ and $\|(a,b)\|_*=1$
            \State $A_{t+1} \gets A_t \cup \{(a,b)\}$.
            \State $G_{t+1} \gets G_t$.
        \EndIf
        \EndFor
        \State \Return $\ub$.
    \end{algorithmic}
\end{algorithm}

In Line~\ref{algdefx}, $\nabla \pot(p_t)[1:n]$ denotes the first $n$ components of $\nabla \pot(p_t)$, and $\nabla \pot(p_t)[n+1]$ denotes the last component of $\nabla \pot(p_t)$.
The sets~$A_t$ and~$G_t$ denote the already known/separated inequalities and objective gradients during iteration~$t$.

\begin{lemma}\label{alwaysANewConstraint}
    When $x_t \in \R^n$ is computed in iteration $t$ of \Cref{algorithm},
    it is well-defined and
    we have $\sprod{c}{x_t} \leq d$ for every $(c,d) \in A_{t} \cup G_t$.
\end{lemma}
\begin{proof}
    Since $p_t$ minimizes $\pot$ over $\conv(A_t \cup G_t)$,
    for every $q \in \conv(A_{t} \cup G_t)$
    we have $\sprod{\nabla \pot(p_t)}{q-p_t} \geq 0$.
    If $p_t \neq \zero$ then from strict convexity of $\pot$ and $\min_{x\in\R^{n+1}}\pot(x)=\pot(\zero)=0$ we get
    \begin{equation}\label{minimizing-f}
        \sprod{\nabla \pot(p_t)}{q} \geq \sprod{\nabla \pot(p_t)}{p_t} > 0.
    \end{equation}
    First, apply~\eqref{minimizing-f} to $q = (\zero,1)/\|(\zero,1)\|_* \in A_t$
    and conclude $\nabla \pot(p_t)[n+1] > 0$. This makes sure that
    $x_t$ can be computed.
    Second, we apply Inequality \eqref{minimizing-f} to $q = (c,d) \in A_t \cup G_t$
    and find that
    \(
        d - \sprod c {x_t} = \frac{1}{\nabla \pot(p_t)[n+1]}\sprod{\nabla \pot(p_t)}{(c,d)} > 0,
    \)
    thus $x_t$ satisfies $\sprod{c}{x_t} \leq d$ for all $(c,d)\in A_t \cup G_t$.
    \qed
\end{proof}

Note that, for the sake of presentation, in Line~\ref{algmin} we require $p_t$ to be the convex combination of minimum $\pot$-value.
However, it is usually not necessary to compute such a minimum.
The same convergence rates can be obtained if, in every iteration,
$p_t$ is a suitable convex combination of $p_{t-1}$ and some
${(c,d) \in A_t \cup G_t}$ with $\sprod{\nabla \pot(p_{t-1})}{(c,d)} < 0$.
If the last coordinate of $p_{t-1}$, as discussed in the above proof,
is not positive, then such an update can be made towards $(\zero,1)/\|(\zero,1)\|_* \in A_t$.
Any such update will significantly decrease $\pot(p_t)$, and the computation
in Line~\ref{algmin} is guaranteed to make at least that much progress.
This shows that simple updates of $p_t$, which may be more preferable in practice, still suffice to achieve the claimed convergence rates.

\begin{lemma}\label{lem:FW}
    Suppose that $\pot$ is $1$-smooth with respect to $\|\cdot\|_*$ and that
    \[
        \|(\nabla f(x), \sprod{\nabla f(x)}{x})\|_* \leq 1
    \]
    for every $x \in K$.
    Then for every $t = 1,\dots,T$, \Cref{algorithm} satisfies $\pot(p_t) \leq \frac{8}{t+2}$.
\end{lemma}
\begin{proof}
    Recall that $\pot$ is $1$-smooth with respect to $\|\cdot\|_*$ if
    \[
        \sprod{\nabla \pot(x)}{x-y} \le \pot(x) - \pot(y) + \tfrac{1}{2} \|x-y\|_*^2
    \]
    holds for all $x,y \in \R^{n+1}$.
    If $p_t \ne \zero$, then since $p_1 = \frac{1}{\|(\zero,1)\|_*} (\zero,1) \in A_t \cup G_t$ we must have $\nabla \pot (p_t)[n+1] > 0$.
    Using this, it is easy to check that in every iteration we add a point $q_t \in A_{t+1} \cup G_{t+1}$ with $\|q_t\|_* \le 1$ such that $\sprod{\nabla \pot(p_t)}{q_t} \le 0$ holds.
    Moreover, note that the $1$-smoothness of $\pot$ yields $ \pot(p_t) \le \pot(p_1) \le \tfrac{1}{2} \|p_1\|_*^2 = \tfrac{1}{2} $.
    Thus, setting $\lambda := \frac{1}{4} \pot(p_t)$ we obtain
    \begin{align*}
        \pot(p_{t+1})
        \le \pot(p_t + \lambda(q_t - p_t))
        & \le \pot(p_t) - \lambda \sprod{\nabla \pot(p_t)}{p_t - q_t} + \tfrac{1}{2} \lambda^2 \|p_t - q_t\|_*^2 \\
        & \le \pot(p_t) - \lambda \sprod{\nabla \pot(p_t)}{p_t} + \tfrac{1}{2} \lambda^2 \|p_t - q_t\|_*^2 \\
        & \le \pot(p_t) - \lambda \pot(p_t) + \tfrac{1}{2} \lambda^2 \|p_t - q_t\|_*^2 \\
        & \le \pot(p_t) - \lambda \pot(p_t) + 2 \lambda^2 \\
        & = \pot(p_t) - \tfrac{1}{8} \pot(p_t)^2,
    \end{align*}
    where the second inequality holds since $\pot$ is $1$-smooth, the fourth inequality follows from convexity since $\pot(\zero) = \zero$, and the last inequality holds since $\|p_t\|_* \le 1$ and $\|q_t\|_* \le 1$.
    From this we can derive $\frac{1}{\pot(p_t)} \ge \frac{1}{\pot(p_1)} + \frac{1}{8}(t-1)$ for all $t$, which yields the claim since $\pot(p_1) \le \frac{1}{2}$. \qed
\end{proof}


The following lemma yields conditions under which a small value of $\pot(p_t)$ implies that $\ub$ is close to the minimum of $f$ over $K$.
Note in particular that it proves that if $\|p_t\| = 0$ then $\ub = \opt$.

\begin{lemma}\label{lem:ConicDistanceVsGap}
    Assume that $\|(x,-1)\| \leq 2$ holds for every $x \in K$,
    and there exist $z \in K$ and $\alpha \in (0,1]$ such that
    $\sprod{(a,b)}{(-z,1)} \geq \alpha \|(-z,1)\| \|(a,b)\|_*$ holds
    for every $(a,b) \in \oracle \cup \{(\zero,1)\}$.
    Moreover, assume that $\|(\nabla f(x), \sprod{\nabla f(x)}{x})\|_* \leq 1$ holds for every $x \in K$.
    If $\|p_T\| \leq \alpha/2$ in
    \Cref{algorithm}, then the returned value
    satisfies $\ub \ge \opt \ge \ub - \frac{4 \|p_T\|_*(1+\alpha)}{\alpha}$.
\end{lemma}
\begin{proof}
    Let $x^* \in K$ minimize $f(x)$ over $x\in K$ and
    let $F \subset [T-1]$ be the set of iterations (except the last one) in which
    $x_t \in K$.
    Now write the point $p_T$ as a convex combination
    \[
        p_T =
        \sum_{(a,b) \in A_T} \lambda_{(a,b)} (a,b)
        + \sum_{t \in F} \gamma_t (\nabla f(x_t), \sprod{\nabla f(x_t)}{x_t})
    \]
    where $\lambda \geq 0, \gamma \geq 0$ and $\|(\lambda,\gamma)\|_1 = 1$.
    Then we have
    \begin{align*}
        &\sum_{t \in F} \gamma_t (f(x_t) -  f(x^*))
         \leq \sum_{t \in F} \gamma_t \sprod{\nabla f(x_t)}{x_t-x^*} \\
        &\quad = \Big\langle\sum_{t \in F} \gamma_t (\nabla f(x_t),\sprod {\nabla f(x_t)}{x_t}),\;(-x^*,1)\Big\rangle \\
        &\quad \leq \Big\langle\sum_{t \in F} \gamma_t (\nabla f(x_t),\sprod {\nabla f(x_t)}{x_t}) + \sum_{(a,b)\in A_T} \lambda_{(a,b)} (a,b),\;(-x^*,1)\Big\rangle \\
        &\quad = \sprod{p_T}{(-x^*,1)} \leq \|p_T\|_* \cdot \|(-x^*,1)\| \leq 2 \|p_T\|_*.
    \end{align*}
    Here, the inequalities respectively arise from convexity of $f$, that $x^* \in K$ satisfies
    $\sprod{(a,b)}{(-x^*,1)} \geq 0$ for every $(a,b) \in A_T$, and the Cauchy--Schwarz inequality.
    In particular, we find that
    $\min_{t \in F} f(x_t) - f(x^*) \leq \frac{2 \|p_T\|_*}{\sum_{t \in F} \gamma_t}$
    whenever ${\sum_{t \in F}\gamma_t > 0}$. To lower bound this latter quantity,
    we use the assumptions on $z$ to derive the inequalities
    \begin{align*}
        & \alpha \left(1-\sum_{t\in F}\gamma_t\right) \|(-z, 1)\| =
        \alpha \|(-z, 1)\| \sum_{(a,b)\in A_T} \lambda_{(a,b)}\\
        & \quad \leq \sprod{\sum_{(a,b)\in A_T} \lambda_{(a,b)} (a,b)}{(-z,1)} \quad \left(\text{ since $\|(a,b)\|_*=1$ }\right) \\
        & \quad = \sprod{p_T}{(-z,1)} - \sum_{t \in F} \gamma_t \sprod{(\nabla f(x_t),\sprod{\nabla f(x_t)}{x_t})}{\;(-z,1)} \\
        &\quad \leq \|p_T\|_*\cdot\|(-z,1)\| + \sum_{t \in F} \gamma_t \|(\nabla f(x_t),\sprod{\nabla f(x_t)}{x_t})\|_* \cdot \|(-z,1)\|.
    \end{align*}
    Now observe that $\|(\nabla f(x_t),\sprod{\nabla f(x_t)}{x_t})\|_* \leq 1$ for every $t \in F$
    and divide through by $\|(-z,1)\|$ to find
    $\alpha (1-\sum_{t\in F}\gamma_t) \leq  \|p_T\|_* + \sum_{t \in F} \gamma_t$.
    Hence, if $\|p_T\|_* \leq \frac \alpha 2 $ then
    $\alpha/2 \leq (\alpha + 1)\sum_{t \in F}\gamma_t$.
    This lower bound on $\sum_{t \in F}\gamma_t$ suffices to prove the lemma.
    \qed
\end{proof}

Combining the previous two lemmas, we obtain the following convergence rate of our algorithm:

\begin{theorem}\label{thm:f}
  Assume that~$\beta > 0$ is such that $\pot(x) \geq \beta \|x\|^2_*$ for all~$x \in \R^{n+1}$.
  Under the assumptions of Lemmas~\ref{lem:FW} and \ref{lem:ConicDistanceVsGap}, \Cref{algorithm} computes, for every $T \geq \frac{32}{\beta \alpha^2}$, a value $\ub < \infty$ satisfying $\ub \geq \min_{x \in K} f(x) \geq \ub - \frac{16}{\sqrt{\beta(T+2)}} \cdot \frac{1+\alpha}{\alpha}$.
\end{theorem}
\begin{proof}
After $T$ iterations, we have $\beta \|p_T\|_*^2 \leq \pot(p_T) \leq
\frac{8}{T+2} \leq \beta \alpha^2/4$ per \Cref{lem:FW}. Since then
$\|p_T\|_* \leq \frac{\sqrt{8}}{\sqrt{\beta(T+2)}} \leq \alpha/2$,
\Cref{lem:ConicDistanceVsGap} tells us that $\opt \geq \ub - \frac{16
(1+\alpha)}{\sqrt{\beta(T+2)}\alpha}$. \qed
\end{proof}

Let us now apply the previous findings to a concrete setting, in which we assume that the objective function $f$ is $L$-Lipschitz, i.e., $|f(x) - f(y)| \le L\|x-y\|_2$ for all $x,y \in \R^n$.

\begin{theorem}
  \label{thmMainConcrete}
Let $K \subset \R^n$ be a convex body satisfying $z + r\mathbb B_2^n \subset
K \subset R\mathbb B_2^n$, given by a separation oracle $\oracle$, and let
$f\colon \R^n \rightarrow \R$ be an $L$-Lipschitz convex function given by a
subgradient oracle.

Apply \Cref{algorithm} to the function $\frac{1}{RL}f$ using norm $\|(x,y)\| \coloneqq \sqrt{2} \|(x/R,y)\|_2$ and potential $\pot(a,b) \coloneqq \frac{1}{4} \|(R a,b)\|_2^2$.
Then, for every $\eps > 0$, after
\[
    T = O \left(\frac{R^2}{r^2} \cdot \frac{R^2L^2}{\eps^2} \right)
\]
iterations we have $\ub \geq \min_{x \in K} f(x) \geq \ub - \eps$.
\end{theorem}
\begin{proof}
By replacing $f(x)$ by $f(Rx)/(RL)$, $K$ by $K/R$, $\eps$ by $\eps/(RL)$, $r$
by $r/R$, and $z$ by $z/R$, we may assume that $R=L=1$, that $r \in (0,1]$. After this
rescaling, note $\|(x,y)\| \coloneqq \sqrt{2} \|(x,y)\|_2$ and $\pot(a,b)
\coloneqq \frac{1}{4} \|(a,b)\|^2_2 = \frac{1}{2} \|(a,b)\|_*^2$. Crucially,
note that \Cref{algorithm} is invariant under the above replacement. 

We now claim that our choice of input satisfies the conditions of
Theorem~\ref{thm:f} with $\beta = 1/2$ and $\alpha = r/4$. Given the claim,
\Cref{thm:f} directly proves the result. To prove the claim, apart from
verifying that the bounds on $\beta$ and $\alpha$ hold, we must verify
smoothness of $\pot$ with respect to the dual norm, a bound of $2$ on the
norm of $(-x,1)$ for $x \in K$, as well as a dual norm bound of $1$ on
$(\nabla f(x), \sprod{\nabla f(x)}{x})$ for $x \in K$. 

The setting $\beta=1/2$ is direct by definition of $\pot$. Since
$\|\cdot\|_*$ is a Euclidean norm, it is immediate that $\pot$ is $1$-smooth
with respect to $\|\cdot\|_*$. For each $x \in K$, using that $R=L=1$, we may
also verify that
\[
\|(x,1)\| = \sqrt{2} \|(x,1)\|_2 = \sqrt{2} \sqrt{\|x\|_2^2+1} \leq \sqrt{2} \sqrt{R^2+1} = 2,
\]
and
\begin{align*}
\|(\nabla f(x),\sprod{\nabla f(x)}{x})\|_* &= \frac{1}{\sqrt{2}}\|(\nabla f(x),\sprod{\nabla f(x)}{x})\|_2 \\ & \leq \frac{1}{\sqrt{2}}\sqrt{\|\nabla f(x)\|_2^2 + \|\nabla f(x)\|^2\|x\|^2} \\ & \leq \frac{1}{\sqrt{2}}\sqrt{L^2 + L^2R^2} = 1.
\end{align*}

We now show the lower bound $\alpha \geq r/4$. Firstly, since
$\|(-z,1)\|\|(\zero,1)\|_* = \|(-z,1)\|_2\|(\zero,1)\|_2 \leq \sqrt{2}$, we
see that $\sprod{(-z,1)}{(\zero,1)} = 1 \geq \frac{1}{2}
\|(-z,1)\|\|(\zero,1)\|_*$. Next, any $(a,b)$ returned by the oracle is
normalized so that $\|(a,b)\|_* = 1 \Leftrightarrow \|(a,b)\|_2=\sqrt{2}$.
Note then that $\|(-z,1)\| \|(a,b)\|_* \leq 2$. From here, we observe that
\[
 \sprod{(a,b)}{(-z,1)} = b - \sprod a z = b - \sprod{a}{z + ra/\|a\|_2} + r\|a\|_2 \geq r\|a\|_2,
\]
since $z+r a/\|a\|_2 \in K$ by assumption. Furthermore, $b - \sprod a z \geq
b - \|a\|_2 \|z\|_2 \geq b - \|a\|_2$ and $0 \leq b - \sprod a z \leq b +
\|a\|_2$. Thus, $b - \sprod a z \geq \max \{ r\|a\|_2, b-\|a\|_2\}$. We now
examine two cases. If $\|a\|_2 \geq 1/2$, then $b - \sprod a z \geq r/2 \geq
r/4 \cdot \|(-z,1)\| \|(a,b)\|_*$. If $\|a\|_2 \leq 1/2$, then $|b| \geq 1$
since $\|(a,b)\|_2^2 = 2$. This
gives $b- \sprod a z \geq b-\|a\|_2 \geq 1/2 \geq r/2$. Thus, $\alpha \geq r/4$, as needed.
    \qed
\end{proof}

\section{Computational experiments}

In this section, we provide a computational comparison of our method with the standard cut loop, the ellipsoid method, and the analytic center cutting plane method on a testbed of linear optimization instances.
For comparison purposes, all four methods are embedded into a common cutting plane framework such that the same termination criteria apply.

\textit{Framework.} Each method has access to a separation oracle that is equipped with a set of initial linear inequalities valid for $K$ (such as bounds on variables), which are incorporated within each method in a straightforward way.
For instance, we initialize our algorithm by adding these constraints to the set $A_1$.
Moreover, for each instance, we will be given a finite upper bound $\ub$ and incorporate the linear inequality $f(x) \le \ub$ in a similar way.
This upper bound gets updated whenever a feasible solution of better objective value was found.
Our framework collects all inequalities queried by the current method and computes the resulting lower bound on the optimum value in every iteration.
Each method is stopped whenever the difference of upper and lower bound is below~$10^{-3}$.

We will also inspect the possibility of a \emph{smart} oracle that, regardless of whether a given point $x$ is feasible, may still provide a valid inequality as well as a feasible solution (for instance, by modifying $x$ in a simple way so that it becomes feasible).
Such an oracle is often automatically available and can have a positive impact on the performance of the considered algorithms.
For the problems we consider, the actual implementation of a smart oracle will be specified below.

\textit{Implementation.}
The framework has been implemented in \code{julia~1.6.2} using \code{JuMP} and \code{Gurobi~9.1.1}.
To guarantee a fair comparison, all four methods have been implemented in a straightforward fashion.
We use the textbook implementation of the ellipsoid method, and Badenbroek's implementation of the analytic center cutting plane method~\cite{badenbroek2020analytic}.
Our method is implemented~\footnote{\url{https://github.com/christopherhojny/supplement_simple-iterative-methods-linopt-convex-sets}} in the spirit of \Cref{thmMainConcrete}, where $p_t$ is computed using \code{Gurobi}.

\textit{Test sets.}
We use three problem classes in our experiments: linear programming formulations of the maximum-cardinality matching problem, semidefinite relaxations of the maximum cut problem, and LPBoost instances for classification problems.

For the maximum-cardinality matching problem, we consider the linear program
\begin{align*}
  \max \ \Big\{ \sum \nolimits_{e \in E} x_e : \ & x \in [0,1]^E, \, \sum \nolimits_{e \in \delta(v)} x_e \le 1 \text{ for all } v \in V, \\
  & \sum \nolimits_{e \in E[U]} x_e \le \tfrac{|U|-1}{2} \text{ for all } U \subseteq V \text{ with } |U| \text{ odd} \Big\},
\end{align*}
due to Edmonds~\cite{Edmonds1965}, where $G=(V,E)$ is a given undirected graph, $\delta(v)$ is the set of all edges incident to~$v$, and~$E[U]$ is the set of all edges with both endpoints in $U$.
The latter constraints are handled within an oracle that computes an inequality minimizing $(|U|-1)/2 - \sum_{e \in E[U]} x_e$, whereas the other inequalities are provided as initial constraints.
For the above problem, the smart version of the oracle does not provide a feasible point since there is no obvious way of transforming a given point into a feasible one.
However, the smart version always provides the minimizing inequality.

We consider 16 random instances with~500 nodes, generated as follows.
For each~$r \in \{30,33,\dots,75\}$ we build an instance by sampling~$r$
triples of nodes~$\{u,v,w\}$ and adding the edges of the induced triangles
to the graph, forming the test set \texttt{matching}.
We believe that these instances are interesting because the~$r$ triangles
give rise to many constraints to be added by the oracle.
Moreover, we selected all~13 instances from the Color02
symposium~\cite{Color02} with less than~300 edges, yielding the test set \texttt{matching02}.

Our second set of instances is based on the semidefinite relaxation of Goemans and Williamson~\cite{GoemansWilliamson1995} for the maximum cut problem
\begin{align*}
  \max \ \Big\{ \sum \nolimits_{\{v,w\} \in E} c(v,w) (1-X_{v,w})/2 : \ & X_{v,w} = X_{w,v} \text{ for all } v,w \in V, \\
  & X_{v,v} = 1 \text{ for all } v \in V, \, X \succeq 0 \Big \}
\end{align*}
where $c$ are edge weights on the edges of $(V,E)$.
We add the box constraints ${X \in [-1,1]^{V \times V}}$ to the initial constraints and handle the semidefiniteness constraint by a separation oracle that, given $X$, computes an eigenvector $h$ of $X$ of minimum eigenvalue and returns the inequality $\langle hh^\intercal, X \rangle \ge 0$.

Within the smart version of the oracle, this constraint is returned regardless of the feasibility of $X$.
If $X$ is not feasible, the semidefinite matrix $\frac{1}{\lambda - 1}X - \frac{\lambda}{\lambda -
  1}I$ is returned, where $\lambda$ denotes the minimum eigenvalue and $I$ the identity matrix.
We generated 10 complete graphs on 10 nodes with edge weights chosen uniformly at random in $[0,1]$.

Our third set of instances arises from LPBoost~\cite{demiriz2002linear}, a classifier algorithm based on column generation.
To solve the pricing problem in column generation, the following linear program is solved:
\[
  \max \Big\{\gamma : (\gamma, \lambda) \in [-1,1] \times [0,D]^n,
  \sprod{\one}{\lambda} = 1,
  \sum_{i=1}^m y_i h(x^i, \omega)\lambda_i \leq -\gamma \text{ for }
  \omega \in \Omega 
  \Big\},
\]
where~$\Omega$ is a set of parameters, for~$i \in [m]$, $x^i$ is a data point
labeled as~$y_i = \pm 1$, $h(\cdot, \omega)$ is a classifier
parameterized by~$\omega \in \Omega$ that predicts the label of~$x^i$
as~$h(x^i, \omega) \in \{-1,+1\}$, and~$D > 0$ is a parameter.
In our experiments, we restrict~$h(\cdot, \omega)$ to be a decision tree of
height~1, so-called tree stumps, and choose~$D = \frac{5}{n}$.
To separate a point~$(\gamma', \lambda')$, we use \code{julia}'s
\code{DecisionTree} module to compute a decision stump with score
function~$\lambda'$ that weights the data points, whose corresponding
inequality classifies~$(\gamma', \lambda')$ as feasible or not.
A smart oracle always returns the computed inequality and
decreases~$\gamma'$ until~$(\gamma', \lambda')$ becomes feasible according
to the found decision stump.

We extracted all data sets from the UC Irvine Machine Learning Repository~\cite{MLRepo} that are labeled as
multivariate, classification, ten-to-hundred attributes,
hundred-to-thousand instances.
Data sets with alpha-numeric values or too many missing values have been
discarded.

\textit{Results.}
In what follows, we report on the number of iterations, i.e., oracle calls, each method needs to obtain a gap (upper bound minus lower bound) below~$10^{-3}$.
We impose a limit of 500 iterations per instance.
Since we are testing naive implementations of each method, we do not report on running time.

To get more insights on the primal and dual performance of the tested
methods, we also report on their \emph{primal and dual integrals}.
Note that we are solving maximization problems in this section, as opposed
to minimization problems in Section~\ref{sec:alg}.
That is, primal (dual) solutions provide lower (upper) bounds on~$\opt$.
If~$\ell_i$ is the lower bound on the optimal objective value~$\opt$ in iteration~$i$, the
\emph{primal integral} is~$\sum_{i = 1}^{500} \frac{\opt - \ell_i}{\opt - \ell_1}$.
The \emph{dual integral} is computed analogously.
If an integral is small, this indicates quick progress in finding the
correct value of the corresponding bound.

\begin{table}[t]
  \begin{scriptsize}
    \caption{Comparison of iterations and dual/primal integral without smart oracles.}
    \label{tab:itercntDInosmartNoinit}
    \begin{tabular*}{\textwidth}{@{}l@{\;\;\extracolsep{\fill}}rrrrrrrrrrr@{}}\toprule
             & \multicolumn{4}{c}{\#iterations} & \multicolumn{4}{c}{dual integral} & \multicolumn{3}{c}{primal integral}\\
      \cmidrule{2-5} \cmidrule{6-9} \cmidrule{10-12}\\
      instance & LP & ellipsoid & analytic & our & LP & ellipsoid & analytic & our & ellipsoid & analytic & our\\
      \midrule
                      matching & 175.44 & 500.00 & 500.00 &  99.81 &  48.34 & 473.02 &  22.13 &  21.10 &  52.12 &   9.29 &   4.40\\
               matching02 & 283.77 & 460.77 & 491.69 &  47.15 & 257.76 & 339.67 & 194.26 &  21.64 &  23.41 &   5.91 &   2.13\\
                        maxcut & 265.30 & 500.00 & 500.00 & 193.30 &   7.72 &  44.32 &   3.48 &   6.14 &  21.15 &   9.04 &   6.32\\
                       LPboost &  91.94 & 489.06 & 479.12 & 278.06 &   3.15 &  13.62 &  20.65 &  53.15 & 459.97 & 100.71 &  64.08\\
      \bottomrule
    \end{tabular*}
  \end{scriptsize}
\end{table}

Table~\ref{tab:itercntDInosmartNoinit} summarizes our results without smart oracles, where all numbers are average values.
Here, ``matching'' refers to the random instances and ``matching02'' to the
instances from the Color02 symposium.
The standard cut loop is referred to as ``LP'', the
ellipsoid method as ``ellipsoid'', the analytic center method as
``analytic'', and \Cref{algorithm} as ``our''.
Note that Table~\ref{tab:itercntDInosmartNoinit} does not report on the
primal integral of ``LP'' since the standard cut loop is a dual
method.

We see that the ellipsoid and analytic center methods are
struggling with solving any instance within~500 iterations independent from
the problem class.
Our algorithm solves the instances of the matching and max-cut
problem much faster than the standard cut loop.
Only for LPBoost, the standard cut loop clearly dominates our algorithm.
To better understand this behavior, the integrals reveal that our algorithm is better in improving the primal bound than the dual
bound, with the only exception being LPBoost.
The analytic center method, however, performs significantly worse than our algorithm in improving the primal bound.
Regarding the dual bound, it performs better than our algorithm (with
the exception of matching02).
The ellipsoid method is much worse in improving the primal bound in
comparison with the analytic center method and our algorithm.
Regarding the dual bound, a similar trend can be observed with LPBoost
being an exception.

\begin{figure}[t]
  \vspace{-1em}
  \centering
  \includegraphics[scale=0.3]{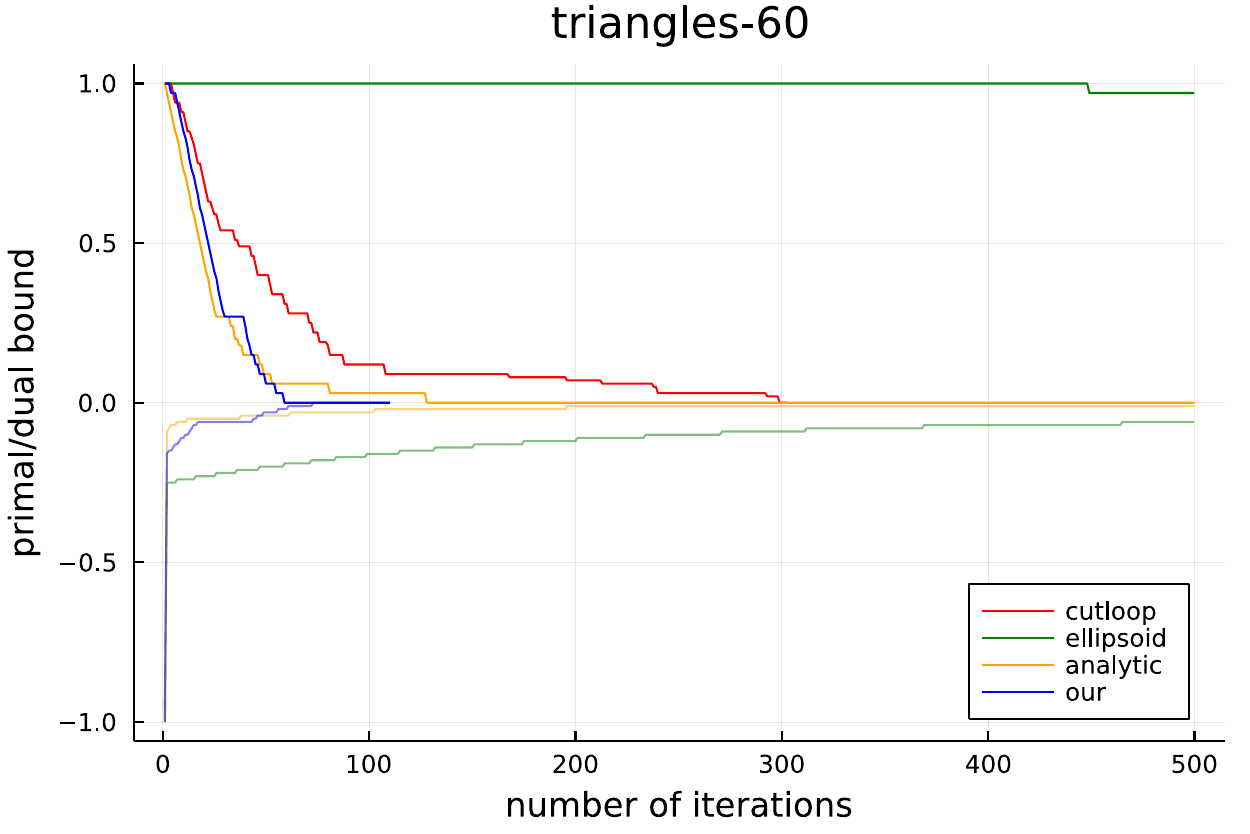}
  \captionof{figure}{Typical primal/dual bounds for a random matching instance.}
  \label{fig:plot}
  \vspace{-1.5em}
\end{figure}

In summary, the analytic center cutting plane method improves the dual
bound more quickly than our algorithm.
It can find a good primal
solution early as the primal integral is small,
however it fails to close the
remaining gap within the iteration limit.
Our algorithm is able to close the primal gap faster, with the trade-off of a slightly slower dual
convergence.
A typical plot of the of the relative primal and dual gaps is given in Figure~\ref{fig:plot}.

\begin{table}[b]
  \vspace{-3em}
  \begin{scriptsize}
    \caption{Comparison of iterations and dual/primal integral with smart
      oracles.}
    \label{tab:itercntDIsmart}
    \begin{tabular*}{\textwidth}{@{}l@{\;\;\extracolsep{\fill}}rrrrrrrrrrr@{}}\toprule
             & \multicolumn{4}{c}{\#iterations} & \multicolumn{4}{c}{dual integral} & \multicolumn{3}{c}{primal integral}\\
      \cmidrule{2-5} \cmidrule{6-9} \cmidrule{10-12}\\
      instance & LP & ellipsoid & analytic & our & LP & ellipsoid & analytic & our & ellipsoid & analytic & our\\
      \midrule
                      matching & 175.44 & 500.00 & 500.00 &  99.81 &  48.34 & 473.02 &  22.13 &  21.10 &  52.12 &   9.29 &   4.40\\
               matching02 & 283.77 & 460.77 & 491.69 &  47.15 & 257.76 & 339.67 & 194.26 &  21.64 &  23.41 &   5.91 &   2.13\\
                        maxcut & 265.30 & 500.00 & 500.00 & 231.00 &   7.72 &  42.90 &   3.48 &   6.15 &  20.42 &   8.91 &   5.59\\
                       LPboost &  86.94 & 346.38 &  88.00 & 127.00 &   3.04 &  13.50 &   5.54 &   5.46 &  25.41 &   6.83 &   6.95\\
      \bottomrule
    \end{tabular*}
    \vspace{-1em}
  \end{scriptsize}
\end{table}

In a second experiment, we investigate the effect of smart oracles.
As Table~\ref{tab:itercntDIsmart} shows, the algorithms mostly benefit from having access to a smart oracle in the case of LPBoost.
A reason might be in the particular structure of these instances: the objective just consists of~$\gamma$ and every truncated convex combination~$\lambda$ is feasible.
However, there is no impact of smart oracles on the matching and maxcut instances, respectively.

\bibliographystyle{splncs04}
\bibliography{references}

\end{document}